# Optimal management of the vaccination process in SIRD epidemic models under constraints

Bogdan Norkin[1] [0009-0001-7531-0420] and Vladimir Norkin[2] [0000-0003-3255-0405]

[1,2] V.M.Glushkov Institute of Cybernetics, Kyiv, Ukraine
bgdan.norkin@gmail.com
vladimir.norkin@gmail.com

**Abstract.** The paper considers the problems of optimal vaccination control in the classical SIR model under constraints on the resource capabilities of the insurance medical system, in particular under constraints on the possible absolute rate of vaccination of the population and the limitation on the available number of vaccines. The application of classical optimal control methods, the dynamic programming method and the Pontryagin maximum principle for such a model encounters difficulties associated with the possible non-smoothness of the Bellman function, and in the Pontryagin method the problem is to solve a boundary value problem with discontinuous control. Therefore, in the paper, optimal control is sought in the class of so-called parametric strategies, which reduces the original problem to a finite-dimensional optimization problem with respect to unknown parameters.

## 1    Introduction

COVID is not receding, the problem of epidemic control and management remains relevant. The classical SIR (Susceptible, Infected, Recovered) model of epidemic development in the form of a system of differential equations was proposed in the work (Kermack, McKendrick, 1927). A review of research on this model and others is in the work (Hethcote, 2000). An introduction to the modern theory of mathematical modeling of epidemics, the nomenclature of deterministic models and methods of optimal epidemic management based on the Pontryagin maximum principle are presented in the books (Kolesin Zhitkova, 2004; Martcheva, 2015). Stochastic models of epidemics are considered in the works (Britton, 2010; Allen, 2008; Andersson, Britton, 2012; Gill, 2015; Britton et al., 2019) (Katriel, 2013; Ishikawa, 2013; Bogdanov, 2023; Atoyev et al., 2024; Bogdanov, Knopov, 2022; Njiasse et al., 2025; Bogdanov, 2022). The use of cellular models to predict the COVID-19 epidemic process in Ukraine is reflected in the works (Brovchenko, 2020; Knopov et al., 2020; Komisarenko, 2020; Nesteruk, 2021; Chumachenko, Chumachenko, 2023; Pantyo et al., 2024). The formulation of optimal control problems of the epidemic process using the Pontryagin maximum principle (Pontryagin, 2018) is considered in the works (Kolesin, Zhitkova, 2004; Gaff, Schaefer, 2009; Laarabi et al., 2013; Ishikawa, 2013;

Martcheva, 2015; Batalin, Terletsky, 2015; Balderrama et al., 2022; Ramponi, Tessitore, 2024; Belili et al., 2024; Karim et al., 2025; Bolzoni, Della Marca, 2025; Wanjala et al., 2025). A review of methods and software for solving boundary value problems in the Pontryagin method is in the work (Mazzia, Settanni, 2021). The application of the dynamic programming method (Bellman, Dreyfus, 2015) to solve this problem is considered in the works (Hethcote, Waltman, 1973; Touffik, Slimani, 2022; Federico et al., 2024; Njiasse et al., 2025), the reduction to a nonlinear programming problem in (Olivares, Staffetti, 2021). The epidemic process management model in the form of a multi-stage stochastic programming model is considered in (Yin, Büyüktahtakın, 2021; Yin, Büyüktahtakın, 2022) and in the form of a stochastic programming model with probabilistic constraints in (Tanner et al., 2008; Gujjula et al., 2022).

The difficulties of applying the dynamic programming method are associated with the possible non-smoothness of the Bellman function, and in the Pontryagin method, the problem is solving a boundary value problem with discontinuous control.

This paper uses a different approach to optimal management of the epidemic process through vaccination, namely, optimal management is sought in some parametric class, which is given by a finite number of unknown parameters. After substituting such parametric management, the original optimal management problem is transformed into a finite optimization problem to find optimal values of unknown parameters. The total costs of the insurance medical system for vaccination and treatment of patients are used as the objective function. Examples of this approach to solving optimal management problems are in the works (Floudas et al., 2013; Knopov, Norkin, 2022). Another innovation of the approach is to take into account the resource constraints of the medical system, restrictions on the speed of vaccination and the availability of the vaccine.

## 2  Classical SIRD model in the form of differential equations

Classic continuous dynamic epidemic SIRD model (Susceptible, Infected, Recovered, Dead) for a population of $N$ individuals includes the following values (Kermack, McKendrick, 1927):

$N$ - total number of individuals; $t \in [0, T]$ - time;

$S(t)$ - the number of healthy individuals (without immunity, Susceptibles);

$I(t)$ - the number of patients and infectious (Infected);

$R(t)$ - the number of individuals with immunity (healthy with immunity and those who recovered with the acquired immunity);

$D(t)$ - the number of dead;

$u(t)$ - the share of vaccinated (isolated) per unit time;

$p(t)dt = p(I(t)) \cdot dt$ - the probability for a healthy person to be infected during time $dt$;

$S(t) \cdot p(I(t))dt$ - average number of infected during time $dt$;

$\alpha \cdot dt$ - the likelihood of recovery of infected during time $dt$;

$\beta \cdot dt$ - the likelihood of death (or serious illness) of infected during time $dt$;

$u(t) \cdot dt$ - (managed) percentage of vaccinated during time $dt$;

$v(t)dt = S(t) \cdot u(t)dt$ - the number of vaccinated during the time $dt$.

**Infection model (probability of infection) in the SIR model** (Katriel, 2013; Gill, 2015). Let

$r$ - the intensity of contacts of an individual, $r \cdot dt$ - the number of contacts during time $dt$;

$\varepsilon$ - the probability of infection when contacting a patient.

Obviously, $I/(N-D) \approx I/N$ - the probability of meeting a patient during any contact. What is the probability of infection during time $dt$, that is, during $r \cdot dt$ contacts?

Note that

$rdt\dfrac{I}{N}$ - the number of contacts with patients during time $dt$,

$(1-\varepsilon)$ - the probability of not getting infected when contacting a patient,

$(1-\varepsilon)^{rdt\frac{I}{N}}$ - the probability of not getting infected during $rdt\dfrac{I}{N}$ contacts,

$1-(1-\varepsilon)^{rdt\frac{I}{N}}$ - the probability of getting infected during $rdt\dfrac{I}{N}$ contacts,

$$1-(1-\varepsilon)^{r\frac{I}{N}dt} = 1-\exp\left(\frac{r\log(1-\varepsilon)I}{N}dt\right) \approx -\frac{r\log(1-\varepsilon)I}{N}dt.$$

Thus, the infection intensity (probability of infection per unit time) is equal to

$$p(I) = -\frac{rI\log(1-\varepsilon)}{N}.$$

**Classical SIRD-model in the form of differential equations** (Kermack McKendrick, 1927). Let us assume that those individuals who recover acquire lasting immunity. Then the entered quantities are related by the following equations:

$S(t) + I(t) + R(t) + D(t) = N$,

$S(t+dt) = S(t) - S(t)p(I(t))dt - v(t)dt$,

$I(t+dt) = I(t) + S(t)p(I(t))dt - \alpha I(t)dt - \beta I(t)dt$,

$R(t+dt) = R(t) + \alpha I(t)dt + v(t)dt$,

$D(t+dt) = D(t) + \beta I(t)dt$.

From here we get:

$$\frac{dS}{dt} = -Sp(I) - v(t),$$

$$\frac{dI}{dt} = Sp(I) - \alpha I - \beta I = Sp(I) - I, \quad \alpha + \beta = 1,$$

$$\frac{dR}{dt} = \alpha I + v(t),$$

$$\frac{dD}{dt} = \beta I, \quad p(I) = -\frac{rI \log(1-\varepsilon)}{N}, \quad v(t) = S u(t), \quad t \in [0,T].$$

Initial conditions for this system: $I(0) > 0$, $S(0) = N - I(0)$.

**Optimal control problem.** Let us denote the following parameters:
  $a$ - the cost of vaccination of one person,
  $b$ - the cost of treatment of one person,
  $c$ - the cost of death (treatment of a severe course of the disease) for one person.

Functionalities to be optimized (minimized):

$$J(u) = a\int_0^T S(t)u(t)dt + (b\alpha + c\beta)\int_0^T I(t)dt \to \min_{u(t) \in U} \quad \text{(the cost of immunization and treatment)},$$

$I(T,u) \to \min_{u(t) \in U}$ (the number of patients at the time $T$);

$$R(T,u) = \int_0^T \beta I(t)dt \quad \text{(the total number of deaths during the time } T \text{.)}$$

The integral criterion of the optimal epidemic control is:

$$J(u(\cdot)) = \underbrace{a\int_0^T S(t)u(t)dt}_{\text{Vaccination costs}} + \underbrace{b\int_0^T \alpha I(t)dt}_{\text{Treatment costs}} + \underbrace{c\int_0^T \beta I(t)dt}_{\text{Treatment costs of seriously ill}} \to \min_{u(\cdot) \in U}$$

**Invariant form (independence from population size) of SIRD model equations.** Let us denote $s = S/N$, $i = I/N$, $\rho = R/N$, $d = D/N$, $p(i) = -ir\log(1-\varepsilon)$, $s(0) = S(0)/N$, $i(0) = I(0)/N$, $\rho(0) = R(0)/N$, $d(0) = D(0)/N$. Then $s(t) + i(t) + \rho(t) + d(t) = 1$.

The model equations take the following form:

$$\frac{ds}{dt} = -sp(i) - v(t),$$

$$\frac{di}{dt} = sp(i) - i,$$

$$\frac{d\rho}{dt} = \alpha i + v(t),$$

$$\frac{d\,d(t)}{dt} = \beta i, \quad p(i) = -ri\log(1-\varepsilon), \quad v(t) = s(t)u(t).$$

The integral criterion of the task based on one person has the form:

$$j(u(\cdot)) = \int_0^T \left( a\ s(t)\ u(t) + \underbrace{(\alpha b + \beta c)}_{d}\ i(t) \right) dt \to \min_{u(\cdot) \in U}$$

$$\underbrace{\phantom{j(u(\cdot)) = \int_0^T \left( a\ s(t)\ u(t) + (\alpha b + \beta c)\ i(t) \right) dt}}_{Total\ costs}$$

In this model, control $u(t)$ can be interpreted either as the proportion of the total number of susceptibles to the infection that is vaccinated per unit of time, or as the probability that a separate healthy person will be vaccinated per unit of time. This model of vaccination management is adopted, for example, in works (Gaff, Schaefer, 2009; Laarabi et al., 2013; Ishikawa, 2013; Martcheva, 2015; Federico et al., 2024; Bolzoni, Della Marca, 2025; Karim et al., 2025; Wanjala et al., 2025; Njiasse et al., 2025).

However, a more realistic model of vaccination is a model that takes into account the resource limitations of the medical system. In it, the number of vaccinated per unit of time, i.e. the vaccination rate is given by the formula:

$$v(t) = \begin{cases} \min\{k, l \cdot s(t)\}, & 0 \leq t \leq \tau, \\ 0, & t > \tau, \end{cases}$$

Where $k$ is the speed capabilities of the medical system carrying out the vaccination, $l$ is the probability that a healthy person agrees and will be vaccinated, $\tau$ is the duration of the vaccination program. Obviously, the duration of the vaccination program should not exceed the (unknown) duration of the epidemic. We will consider the quantities $k$ and $l$ to be fixed and known. Another ingredient of the model is the restriction on the availability of the vaccine, $\int_0^\tau v(t) dt \leq m$, where $m$ is the available amount of the vaccine.

Thus, the optimization model has the form:

$$j(\tau) = \int_0^\tau \left( a\ v(t) + (b + c\beta)\ i(t) \right) dt \to \min_\tau, \quad (1)$$

where $\tau \leq T$,

$$p(i) = -ri \log(1 - \varepsilon),\ v(t) = \begin{cases} \min\{k, l \cdot s(t)\}, & 0 \leq t \leq \tau, \\ 0, & t > \tau, \end{cases} \quad (2)$$

$$\frac{ds}{dt} = -sp(i) - v(t), \quad (3)$$

$$\frac{di}{dt} = sp(i) - i, \quad (4)$$

$$\frac{d\rho}{dt} = \alpha i + v(t), \quad (5)$$

$$\frac{d\,d(t)}{dt} = \beta i ,\quad (6)$$

$$\int_0^\tau v(t)dt \leq m .\quad (7)$$

This problem is a nonlinear one-dimensional conditional optimization problem with a difficult to calculate criterion. $\tau$ is an unknown quantity, we will denote $\tau^*$ its optimal value. We emphasize that this optimal value was obtained taking into account the total costs of vaccination and treatment. Having solved this problem $\tau^* \leq T$, it is possible to predict the course of the epidemic for the entire planning period $[0,T]$ with limited resources of the medical system, taking into account the total costs of vaccination and treatment. At the same time, it is possible to calculate all the accompanying indicators of the epidemic process: the value of the peak of the epidemic $\max_{0 \leq t \leq T} i(t)$, the time $\arg\max_{0 \leq t \leq T} i(t)$ of the onset of the peak of the epidemic, the duration $t^* = \{t : i(t) \approx 0\} = \tau^*$ of the epidemic, financial costs $j(\tau^*)$, etc., and then make decisions on changing resource parameters $k$ and $m$.

Another possible application of the considered model can be in the planning of vaccine purchases. If we solve model (1)-(6), i.e. without resource limitation (7), then we get a forecast picture of the epidemic and the duration of the epidemic (vaccination period $\tau^{**}$) only if the medical system is limited by the rate of vaccination $k$ and the probabilities $l$ of vaccination propensity. This will allow you to calculate the necessary amount of vaccine $m^{**} = \int_0^{\tau^{**}} \min\{k, s(t)\} dt$ for the implementation of the scenario with a limitation only on the rate of vaccination, as well as calculate all the accompanying indicators of the epidemic process for subsequent decision-making on changing the model parameters.

## 3   Numerical experiments

Numerical experiments to solve problems (1)-(7) were carried out in the Matlab system (Davis, Sigmon, 2005) using functions ode23 (solution of ordinary differential equations) and fminbnd (one-dimensional optimization) with default algorithmic parameters. The following values of model parameters from table 1 were used in numerous experiments. Variants 1 and 2 differ only in the value of one parameter.

Table 1. Numerical values of parameters in model (1)-(7).

| parameters | Variant 1 | Variant 2 | Meaning of a parameter |
|---|---|---|---|
| $\alpha$ | 0.95 | 0.95 | Probability of recovery of an infected person per unit of time |
| $\beta$ | 0.05 | 0.05 | Probability of death (or severe illness) of an infected person per unit of time |
| $r$ | 10 | 10 | Individual contact intensity |
| $\varepsilon$ | 0.3 | 0.3 | Probability of infection upon contact with a sick person |
| $s(0)$ | 0.999 | 0.999 | Initial proportion of susceptible to infection |
| $i(0)$ | 0.001 | 0.001 | Initial proportion of infected |
| $\rho(0)$ | 0 | 0 | Initial proportion of non-susceptible to infection |
| $d(0)$ | 0 | 0 | Initial proportion of deceased |
| $T$ | 15 | 15 | Modelling time period |
| $a$ | 5 | 5 | Cost of one dose of vaccine |
| $b$ | 50 | 50 | Cost of treatment of a mild patient |
| $c$ | 500 | 500 | Cost of treatment of a severe patient |
| $k$ | 0.1 | 0.1 | Vaccination rate |
| $l$ | 0.3 | 0.3 | Proportion of consent to vaccination |
| $m$ | 2.949 | 0.5 | Number of vaccine doses available |

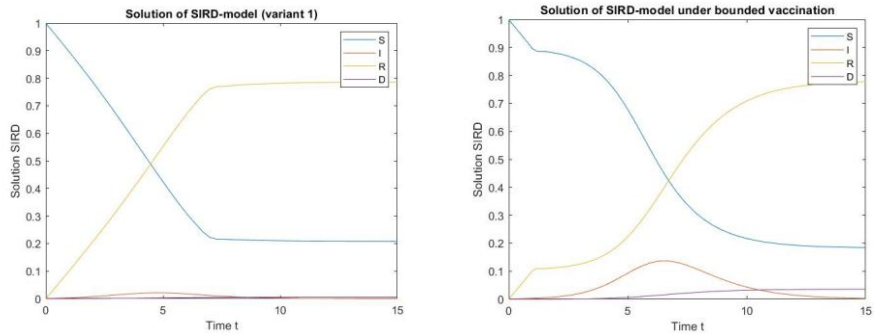

**Fig. 1.** Simulation of the SIRD model functions $s(t), i(t), \rho(t), d(t)$ for variants 1 and 2.

Fig. 1 shows the trajectories of the model variables $s(t), i(t), \rho(t), d(t)$ as a function of time for parameter values from Table 1. The left part of the figure corresponds to the value $m \geq 2.9491$ (option 1, the limit on the number of vaccines is actually absent). The right part of this figure corresponds to the value $m = 0.5$ (option 2, restriction (7) on the number of vaccines is active).

Fig. 2 displays optimal controls as a function of time for options 1 and 2 of model parameter values. As can be seen from the figures, the behavior of the functions $s(t), i(t), \rho(t), d(t)$ of the models and the optimal controls $v(t)$ are radically different in the absence and presence of limitation (7) on the number of available vaccines.

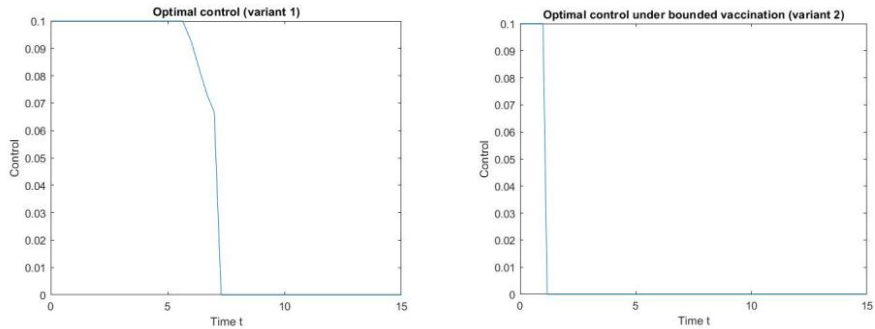

**Fig.** 2. Optimal controls $v(t)$ in variants 1 and 2.

## 4  Conclusions

The paper considers the problem of optimal control of the vaccination process in the SIRD-model under resource limitations of the medical system, namely, when limited to the absolute rate of vaccination and limited to the available amount of vaccines. The criterion of optimality is the total costs of the medical insurance system for vaccination of healthy and treatment of sick individuals. This model differs from the majority of models considered in the literature, in which the limitation on the rate of vaccination is a certain proportion of the number of healthy individuals, which may exceed the actual capabilities of the medical system. As a rule, the limit on the available amount of vaccines is not taken into account. Such simplifying assumptions make it easier to obtain optimal control conditions, but they may be far from reality. Another element of the novelty of this article's approach is the search for optimal control in a certain class of parametric strategies. As a result, the problem of optimal control is reduced to the problem of finite-dimensional optimization with respect to unknown parameters, which define a class of parametric strategies.